\documentclass[a4]{amsart}

\title{Ambient metric construction of CR invariant differential operators}
\author{Kengo Hirachi}
\thanks{Graduate School of Mathematical Sciences, The University of Tokyo, 3-8-1 Komaba, Megro, Tokyo 153-8914 JAPAN}

\newtheorem{theorem}{Theorem}
\newtheorem{proposition}{Proposition}
\newtheorem{lemma}{Lemma}

\def\nelongarrow{\begin{picture}(15,10)(0,0)
\put(2,0){\vector(3,1){20}}
\end{picture}}
\def\selongarrow{\begin{picture}(15,10)(0,0)
\put(2,8){\vector(3,-1){20}}
\end{picture}}

\newcommand{\C}{\mathbb{C}}

\newcommand{\Z}{\mathbb{Z}}

\newcommand{\T}{\mathbb{T}}

\newcommand{\calH}{\mathcal{H}}

\DeclareMathOperator{\image}{image}

\renewcommand{\Re}{\Real}

\renewcommand{\Im}{\Imag}

\newcommand{\MA}{Monge-Amp\`ere }
\renewcommand{\c}{\overline}

\newcommand{\pa}{\partial}
\newcommand{\wt}{\widetilde}

\renewcommand{\th}{\theta}

\newcommand{\bC}{\mathbb{C}}
\newcommand{\bR}{\mathbb{R}}
\newcommand{\bZ}{\mathbb{Z}}

\newcommand{\contr}{\operatorname{contr}}
\renewcommand{\Re}{\operatorname{Re}}
\renewcommand{\Im}{\operatorname{Im}}

\newcommand{\calO}{\mathcal{O}}
\newcommand{\calI}{\mathcal{I}}
\newcommand{\calJ}{\mathcal{J}}
\newcommand{\calE}{\mathcal{E}}
\newcommand{\calN}{\mathcal{N}}
\newcommand{\calP}{\mathcal{P}}

\newcommand{\conj}{\overline}

\begin{document}

\maketitle
\vspace{-10pt}

\begin{center}{\em Dedicated to the memory of Tom Branson
}
\end{center}

\section{Introduction}
These notes are based on my lectures at IMA, in which I tried to explain basic ideas of the ambient metric construction by studying the Szeg\"o kernel of the sphere.
The ambient metric was introduced in Fefferman \cite{F} in his program of describing the boundary asymptotic expansion of the Bergman kernel of strictly pseudoconvex domain. This can be
seen as an analogy of the description of the heat kernel asymptotic in terms of
local Riemannian invariants. The counterpart of the Riemannian invariants for the
Bergman kernel is invariants of the CR structure of the boundary.
Thus the program consists
of two parts:

\smallskip

(1) Construct local invariants of CR structures;

(2) Prove that (1) gives all invariants by using the invariant theory.

\smallskip
\noindent
In the case of the Szeg\"o kernel, (1) is replaced by the construction of local invariants of the Levi form that are invariant under scaling by CR pluriharmonic functions.
We formulate the class of invariants in Sections 2 and 3.
To simplify the presentation, we confine ourself to the case
of the sphere in $\C^n$. It is the model case of 
the ambient metric construction and the basic tools already appears in this setting. We construct invariants (formulated as CR invariant differential operators) by using the ambient space in Section 4
and then explain, in Section 5, how to prove that we have got all.

The CR invariant operators studied in these notes look similar to 
the invariants of CR densities 
studied in \cite{EG1} and \cite{GG}.
However, there is a crucial difference. We here deal with
operators acting on a quotient of
the space of CR densities. They naturally
arise in the description of the jets of geometric 
structures. For example, the jets of CR structures are described by Moser's normal form, which is defined as a slice of the quotient space of CR densities.
In section 6,  we explain that the quotient space
can be realized as cohomology of a subcomplex of Bernstein-Gelfand-Gelfand (BGG) complex. In particular, Moser's normal form
corresponds to the cohomology of the deformation
complex, a subcomplex of the BGG complex for the adjoint representation. 
The space we use for the Szeg\"o kernel comes from the
Rumin complex, the BGG complex for the trivial representation.

General theory of the ambient metric construction in CR case 
has been 
developed in \cite{F}, \cite{BEG} and \cite{H2};
see \cite{H3} for a survey and more comprehensive references.
The proofs of the results in these notes (without citations)
are given in my forthcoming paper, which will also deal with the
general strictly pseudoconvex manifolds. See also \cite{H1}
for the case of 3-dimensions.

\section{Transformation rule of the Szeg\"o kernel}
To motivate the problem, we first recall the definition of the
Szeg\"o kernel and its transformation rule under scaling of
the volume form. (One can skip this section and jump to the formulation of the problem in the next section.)

Let $\Omega\subset\bC^n$ be a bounded domain with smooth boundary $\pa\Omega$ and $\calO(\Omega)$ be the space of holomorphic functions on $\Omega$. We fix a volume form $d\sigma$ on $\pa\Omega$ and define a pre-Hilbert space
$
CH(\Omega)=C^\infty(\c\Omega)\cap\calO(\Omega)
$ 
with inner product $(f_1,f_2)=\int_{\pa\Omega}f_1\c{f_2}d\sigma$.
The completion of $CH(\Omega)$ is called the Hardy space and denoted by $H^2(\Omega,d\sigma)$. As a set, $H^2(\Omega,d\sigma)$
is independent of the choice of $d\sigma$
but the inner product does depend on the choice.
Each element of $H^2(\Omega,d\sigma)$ is realized as
a holomorphic function on $\Omega$ that has $L^2$ boundary value;
thus we can identify $H^2(\Omega,d\sigma)$ with a subspace of
$\calO(\Omega)$.

Choose a complete orthonormal system $\{\varphi_j(z)\}_{j=1}^\infty$ of $H^2(\Omega,d\sigma)$ and make a series
$$
 K_{d\sigma}(z)=\sum_{j=1}^\infty |\varphi_j(z)|^2.
$$
It converges for $z\in\Omega$ and define a smooth function
on $\Omega$, which is independent of the choice of 
$\{\varphi_j(z)\}_{j=1}^\infty$.  The function 
$ K_{d\sigma}(z)$ is called the Szeg\"o kernel of
 $H^2(\Omega,d\sigma)$. 

For the unit ball
$\Omega_0=\{z\in\bC^n:|z|^2=|z_1|^2+\cdots+|z_n|^2<1\}$
with the standard ($U(n)$-invariant) volume form $d\sigma_0$
on the boundary $S^m$, $m={2n-1}$, we can compute the Szeg\"o kernel
by choosing the monomials of $z_1,z_2,\dots,z_n$
as a complete orthogonal system. Summing up $|z^\alpha|^2$ after normalization gives
$$
 K_{d\sigma_0}(z)=c_n (1-|z|^2)^{-n},
$$
where $c_n$ is the inverse of the volume of the sphere $S^m$.
This computation utilizes the symmetry of $S^m$ and
can be applied only to very limited cases. 
Even if one just perturbs the volume form on $S^m$,
it becomes very hard to compute the Szeg\"o kernel.

For small perturbation of the ball and the volume form, 
we can give the form of the boundary asymptotics of the Szeg\"o kernel.

\medskip

\begin{theorem}[Fefferman, Boutet de Monvel--Sj\"ostrand]\label{FBS}
If $\Omega$ is strictly pseudoconvex domain with a smooth defining function $\rho$, $\Omega=\{\rho>0\}$,
and $d\sigma$ be a smooth volume form on $\pa\Omega$, then 
$$
K_{d\sigma}(z)=\varphi(z)\rho(z)^{-n}+\psi(z)\log\rho(z),
\quad
\varphi,\psi\in C^\infty(\c\Omega).
$$
\end{theorem}

The proof is based on the theory of singular integral operators or Fourier integral operators, and it is not practical to use these calculus to compute the expansion explicitly.
We here try to apply the invariant theory to write down the
coefficients $\varphi$ and $\psi$ in terms
of the ``curvature'' of $(\pa\Omega,d\sigma)$.
For this purpose, we derive the transformation rule
of the Szeg\"o kernel under biholomorphic maps and the scaling of
the volume form.

Let $F\colon\Omega\to\Omega'$ be a biholomorphic map between
strictly pseudoconvex domains with smooth boundaries. Then $F$ can be
extended smoothly to the boundary by Fefferman's theorem. 
Suppose that we are given volume forms $d\sigma$ on $\pa\Omega$ and
$d\sigma'$ on $\pa\Omega'$ such that $F^*(d\sigma')=e^fd\sigma$
for a CR-pluriharmonic function $f$ (i.e., $f$ can be extended to
a pluriharmonic function $\wt f$ on $\Omega$).
If $K_{d\sigma}$ and $K_{d\sigma'}$ are respectively the Szeg\"o kernels of
$(\Omega,d\sigma)$ and $(\Omega',d\sigma')$, then
$$
K_{d\sigma'}(F(z))=e^{-\wt f(z)}K_{d\sigma}(z),\quad z\in\Omega.
$$
In particular, the coefficients of the logarithmic terms
$\psi$ and $\psi'$ of $K_{d\sigma}$ and $K_{d\sigma'}$
satisfy
\begin{equation*}\label{trans-rule-psi}
\psi'(F(z))=e^{-\wt f(z)}\psi(z)\mod O(\rho^\infty).
\end{equation*}
On the boundary, we may multiply each side by the
corresponding volume form and eliminate the scaling factor:
\begin{equation}\label{trans-rule-psi}
F^*(\psi'd\sigma')=\psi d\sigma \quad\text{on }\pa\Omega.
\end{equation}
Our aim here is to classify the differential forms that satisfy 
this transformation rule. (We also want to study higher order term
in $\psi$ and $\varphi$, but it requires more geometric tools;
see \cite{H2}.)

\section{Formulation of the problem}
We now restrict the domain to the ball $\Omega_0$ and express 
the transformation rule \eqref{trans-rule-psi} in terms of the
automorphism group of $\Omega_0$.

The automorphisms of $\Omega_0$ are given by the action of $G=SU(1,n)$ as rational maps.
The action is defined by using the embedding of
$\Omega_0$ into the projective space
$$
\Omega_0\ni z\mapsto [1:z]=[\zeta_0:\zeta_1:\cdots:\zeta_n]\in \mathbb{CP}^n.
$$
The ball is then defined by 
$$
L(\zeta)=|\zeta_0|^2-|\zeta_1|^2-\cdots-|\zeta_n|^2>0
$$
and the complex linear transformations that preserve the hermitian form
$L$ give automorphisms of $\Omega_0$.
Let $\wt e_0={}^t(1,1,0,\dots,0)\in\bC^{n+1}$ and
$P=\{h\in G: h\,\wt e_0=\lambda \wt e_0\}$.
Then the action of $P$ fixes $e_0=(1,0,\dots,0)\in S^m$ and we may
write $S^m=G/P$. Note that $P$ is a parabolic subgroup of $G$ and the representation theory of $P$ plays an essential role in the characterization of invariant differential operators.
For $p,q\in\Z$, we denote by
$\sigma_{p,q}$ the (complex) one-dimensional representation
of $P$, where $h\in P$ is represented by $\lambda^{-p}\conj{\lambda}{}^{-q}$. In case $p=q$, we regard $\sigma_{p,p}$ as
a real representation.

The action of $G$ on $\bC^{n+1}\setminus\{0\}$ can be also
seen as isometries of the  Lorentz-hermitian metric
$$
\wt g=d\zeta_0d\c\zeta_0-d\zeta_1d\c\zeta_1-
\cdots-d\zeta_nd\c\zeta_n.
$$
The Lorentz-hermitian manifold $(\bC^{n+1}\setminus\{0\},\wt g)$ is called {\em the ambient space} of $S^m$. The ambient space
admits $\C^*$ action $\zeta\mapsto\lambda\zeta$ and
we may define homogeneous functions on it:
we say that $\wt f(\zeta)\in C^\infty(\bC^{n+1}\setminus\{0\},\bC)$ is homogeneous of degree
$(p,q)\in\bZ^2$ if $\wt f(\lambda \zeta)=\lambda^p\conj{\lambda}{}^q$
for any $\lambda\in\C^*$.
 We can also define homogeneous functions on the real hypersurface $\calN=\{L(\zeta)=0\}\subset\bC^{n+1}\setminus\{0\}$. Let
$$\calE(p)
=\{f(\zeta)\in C^\infty(\calN,\bR):f(\lambda \zeta)=|\lambda|^{2p}f(\zeta)
,\ \lambda\in\C^*\},
 $$
which is a $G$-submodule of $C^\infty(\calN,\bR)$.
Since $\calN$ is  
a $\bC^*$-bundle over $S^m=G/P$, we can identify 
$\calE(p)$ with the space of the sections of
the real line bundle induced from the representation $\sigma_{p,p}$.
In particular,  $\calE(0)$ can be identified with $C^\infty(S^m,\bR)$,
and $\calE(-n)$ is identified with the space of volume forms $C^\infty(S^m,\wedge^m S^m)$. The correspondence of the latter is given by
$$
f(\zeta_0,\zeta_1,\dots,\zeta_n)\mapsto f(1,z)d\sigma_0.
$$
To simplify the notation we write $\calE=\calE(0)$.
Since the space of CR pluriharmonic functions is preserved by
the $G$-action, we may define the submodule 
$$
\calP=\{f\in \calE: \text{$f$  is CR pluriharmonic}\}.
$$ 

When we study the jets of $\calE(p)$ at $e_0$, it is useful to have  coordinates centred at $e_0$. Let
$$
 \xi_0=\zeta_0+\zeta_n,\ \xi_1=\sqrt{2}\zeta_1,
 \dots, \xi_{n-1}=\sqrt{2}\zeta_{n-1},\xi_n=
 \zeta_0-\zeta_n.
$$
Then in the coordinates $(\xi_0,\dots,\xi_n)$, we have $\wt e_0={}^t(1,0,\dots, 0)$
and 
$$
2L(\xi)=\xi_0\conj\xi_n+\xi_n\conj\xi_0-\sum_{j=1}^{n-1}|\xi_j|^2.
$$ 
Hence setting $w_j=\xi_j/\xi_0,\ j=1,\dots,n$, we may realize $\Omega_0$ as the Siegel domain
$
2\Re w_n>\sum_{j=1}^{n-1}|w_j|^2
$
with $e_0=0$.

Now we return to the Szeg\"o kernel.  On the sphere, each volume form is
written as $d\sigma =e^f d\sigma_0$ for $f\in\calE$. Hence the Szeg\"o kernel
can be seen as a functional of $f$. Taking the defining function
$1-|z|^2$, we write
$$
K_{d\sigma}=\varphi_f (1-{z}^2)^{-n}+\psi_f\log (1-|z|^2).
$$
Here we put the subscript $f$ on $\varphi$ and $\psi$ to emphasis that these coefficients are determined by $f$.
 In particular, we can define a map
$$
\Psi\colon \calE\to\calE(-n),\quad
\Psi(f)=\psi_f d\sigma.
$$
The proof of Theorem \ref{FBS} also implies that $\Psi$ is a differential operator, which polynomially depends on the jets of $f$.
The transformation rule \eqref{trans-rule-psi} can be now reformulated as follows:

\medskip

(i) $\Psi$ is $G$-equivariant;

\medskip

(ii) $\Psi(f+h)=\Psi(f)$ for any $h\in\calP$.

\medskip

\noindent
The property (i) follows from the fact that any automorphism $F$
satisfies $F^*(d\sigma_0)=e^f d\sigma_0$ for an $f\in\calP$; (ii)
is the special case of \eqref{trans-rule-psi} in which $F$ is the identity map.
Since the sphere is the model of CR manifolds,
we will use the terminology ``CR invariant" instead of
$G$-equivariant. 
Then our problem can be stated as follows.

\medskip

{\sc Problem 1.} Write down all CR invariant differential operators 
from $\calE$ to $\calE(-n)$ that are invariant under the additions of $\calP$.

\medskip

To be more precise, we should assume that the operator polynomially depends on the jets of $f$, that is, the value of the
operator at $e_0$ is given by a polynomial in the Taylor coefficients
of $f$ at $w=0$. We will always assume this property
in the rest of these notes.

\section{Construction of CR invariant differential operators}

We first study linear CR invariant differential operators.
Such operators define homomorphisms of $G$-modules
and we may apply the results of the representation theory.
Taking the jets at $e_0$, we may reduce the
study of invariant differential operators to
the one for the homomorphisms between generalized Verma modules.
These homomorphisms can be classified by using the affine action of the Weyl group on the highest weight. In particular,  
we obtain the following 

\medskip

\begin{theorem} \label{linear-unique}
There is a unique, up to a constant multiple,  linear CR invariant  differential operator  $\calE\to\calE(-n)$.
\end{theorem}

\medskip

This theorem ensures that we have a unique linear
CR invariant differential operator, but does not provide
its explicit formula.
We construct the linear operator
by using the Laplacian for the ambient metric
$$
\Delta=\pa_{\zeta_0}\pa_{\conj \zeta_0}-\pa_{\zeta_1} \pa_{\conj \zeta_1}-\cdots-\pa_{\zeta_n} \pa_{\conj \zeta_n}.
$$
This is an analogy of the construction of conformally invariant differential operators in \cite{GJMS}.
\begin{theorem} For $f\in\calE$, take
an $\wt f\in C^\infty(\bC^{n+1}\setminus\{0\})$ homogeneous of degree $(0,0)$ such that $\wt f|_\calN=f$.
Then 
$$
(\Delta^n\wt f)|_{\calN}
$$
depends only on $f$ and defines a differential operator 
$\calE\to\calE(-n)$.
\end{theorem}

\medskip

\begin{proof} 
Since $\Delta^n$ is linear, it suffices to show that
 $\wt f|_{\calN}=0$ implies $(\Delta^n\wt f)|_{\calN}=0$. Write $\wt f=Lh$ with a function $h$ of homogeneous degree $(-1,-1)$.
Then we have, for a positive integer $k$, that
$$
\begin{aligned}
\Delta^k(Lh)&=[\Delta^k,L]h+L\Delta^k h\\
&=k(Z+\c Z+n+k)\Delta^{k-1}h+O(L)\\
&=
k(n-k)\Delta^{k-1}h+O(L).
\end{aligned}
$$
Here $Z=\sum_{j=0}^n \zeta_j\pa_{\zeta_j}$ and
we have used the fact that $\Delta^{k-1}h$
is homogeneous of degree $(-k,-k)$.
In particular, substituting $k=n$, we get $\Delta^n(Lh)=O(L)$.
\end{proof}

\medskip

We denote the linear CR invariant differential operator defined above
by $Q(f)$, as it agrees with the 
CR $Q$-curvature of the volume form 
 $d\sigma=e^f d\sigma_0$ or the contact form $\th$ satisfying 
$d\sigma=\th\wedge (d\th)^{n-1}$. See [FH] for more discussion 
about the ambient metric construction of the $Q$-curvature. 
Since $h\in\calP$ can be extended to 
a pluriharmonic function on the ambient space,
we have $Q(h)=0$ for $h\in\calP$; thus $Q$ satisfies the condition of Problem 1.

We can also use the ambient Laplacian $\Delta$ to construct non-linear CR invariant differential operators. In this case, we need
to specify the ambient extension $\wt f$ more precisely. 

\medskip

\begin{lemma}
Each $f\in\calE$ can be extended to a smooth homogeneous function
$\wt f$ on $\C^{n+1}\setminus \{0\}$ satisfying
\begin{equation}\label{harmonic}
\Delta\wt f=O(L^{n-1}).
\end{equation}
Such an $\wt f$ is unique modulo $O(L^n)$.
\end{lemma}

\medskip

\begin{proof}
We construct $\wt f$ by induction. Suppose that
we have a function $\wt f_k$ homogeneous of degree $(0,0)$
such that $\wt f_k=f+O(L)$ and $\Delta\wt f_k=L^{k-1}h$
for an $h$ homogeneous of degree $(-k,-k)$. (The second condition is vacuous if $k=1$ and we have $\wt f_1$.) We set
$\wt f_{k+1}=\wt f_k+ L^{k}g$ for 
$g$  homogeneous of degree $(-k,-k)$ and try to achieve $\Delta \wt f_{k+1}=O(L^{k})$. Since
$$
[\Delta,L^k]=k\,L^{k-1}(Z+\c Z+n+k), 
$$
we have
$$
\begin{aligned}
\Delta \wt f_{k+1}&=\Delta \wt f_k+[\Delta, L^{k}]g+O(L^{k})\\
&=L^{k-1}h+k(n-k)L^{k-1}g+O(L^{k}).
\end{aligned}
$$
Thus, as long as $k\ne n$, we may set $g=-h/(k(n-k))$ and get $\wt f_{k+1}$. The inductive step starting with $\wt f_1$ gives
$\wt f_{n}$. The uniqueness of $\wt f$ is also clear from the construction.
\end{proof}

\medskip

From this proof, one can also observe that $Q(f)$ is the only obstruction to the existence of an exact formal solution to $\Delta\wt f=0$. In fact, applying $\Delta^{n-1}$ to  $\Delta\wt f=L^{n-1}h$, one gets $\Delta^n\wt f=\text{const.}\, h+O(L)$. Thus
$Q(f)=0$ implies $h=O(L)$ and so $\Delta\wt f=O(L^n)$.
We can then apply the inductive step to get $\wt f_k$
for any $k>n$.

The metric $\wt g$ defines a flat hermitian connection
$\wt \nabla=\wt \nabla^{1,0}+\wt \nabla^{0,1}$. For iterated derivatives $\wt \nabla^{p+q}$, we denote by $\wt \nabla^{p,q}$ the part of type $(p,q)$.  Then we can define scalar valued nonlinear differential operators for $\wt f$ by making complete contractions:
\begin{equation}\label{Weyl-def}
W(\wt f)=\contr( \wt \nabla^{p_1,q_1}\wt f\otimes\cdots\otimes 
\wt \nabla^{p_r,q_r}\wt f)|_{\calN}.
\end{equation}
Here $r\ge2$, $p_j,q_j\ge1$, $p_1+\cdots+p_r=q_1+\cdots +q_r=n$ and the contraction is taken with respect to $\wt g$ for some paring of holomorphic and anti-holomorphic indices.
In view of the ambiguity of $\wt f$, we can show that
$W(\wt f)$ depends only on $f$ and define a CR invariant differential operator
$W\colon\calE\to\calE(-n)$.
Since $\wt\nabla^{1,1}$ kills pluriharmonic functions, we also see that
$W(f+h)=W(f)$ for any $h\in\calP$.
We now give a solution to Problem 1.

\medskip

\begin{theorem}\label{P-invariants-thm}
Let
$
S:\calE\to 
\calE(-n)
$ be a CR invariant differential operator that satisfies
$S(f+h)=S(f)$ for $h\in\calP$.
Then $S$ is a linear combination of $Q(f)$ and 
complete contractions of the form \eqref{Weyl-def}.
\end{theorem}

\medskip

In case $n=2$, the theorem implies
$$
\Psi(f)=c_2 \,Q(f)+c_2'\,\| \wt \nabla^{1,1}\wt f\|^2,
$$
where $\| \wt \nabla^{p,q}\wt f\|^2=\contr(\wt \nabla^{p,q}\wt f\otimes \wt \nabla^{q,p}\wt f)$.
However, we will see in the next section that $\| \wt \nabla^{1,1}\wt f\|^2=0$ by using the fact that the rank of $T^{1,0}S^m$ is one.
Hence $\Psi(f)=c_2 Q(f)$. The constant $c_2=1/(24\pi^2)$ is identified in \cite{H1}; see also \cite{FH}.

In case $n=3$, we have
$$
\Psi(f)=c_3\, Q(f)+c_3'\,\| \wt \nabla^{2,1}\wt f\|^2+
c_3''\,\contr(\wt \nabla^{1,1}\wt f\otimes \wt \nabla^{1,1}\wt f\otimes \wt \nabla^{1,1}\wt f),
$$
where each term is nontrivial. We know $c_3\neq0$ but haven't
computed other constants.

\section{Jet isomorphism theorem and invariant theory}

We  explain how we reduce Theorem \ref{P-invariants-thm}
to a purely algebraic theorem of representation theory.
Let $J\calE$ denote the space of $\infty$-jet of smooth functions
at $e_0\in S^m$. If one fixes coordinates
of $S^m$ around $e_0$, $J\calE$ is identified with the space
of formal power series centered at $e_0$.
We denote by $J\calP$ the subspace of $J\calE$ consisting of
jets of CR-pluriharmonic functions in a neighborhood of $e_0$.

It is clear that a CR invariant (or $G$-equivariant) differential operator
defines a $P$-equivariant map
$$
I\colon J\calE/J\calP\to\sigma_{-n,-n}.
$$
Conversely, by Frobenius reciprocity, we can extends $I$ uniquely
to a $(\mathfrak{g},P)$-equivariant map
$
\wt I\colon J\calE/J\calP\to J\calE(-n)
$
and it extends to a $G$-equivariant differential operator
on $\calE/\calP\to\calE(-n)$. 
Thus Problem 1 is reduced to

\medskip

{\sc Problem 2.}
{\em  Write down all $P$-equivariant map
$J\calE/J\calP\to\sigma_{-n,-n}.
$ 
}
\medskip

Let us call ``$P$-equivariant map" simply
``CR invariant."
When we study the CR invariants of the quotient module
$J\calE/J\calP$, it is useful to have a slice of the coset.
We here follow Moser's argument \cite{CM} of making normal form for
real hypersurfaces in $\C^n$.

In the coordinates $(w',w_n)\in\C^{n-1}\times\C$, the sphere
$S^m$ is given by $2\Re w_n=|w'|^2$ and $(w',v)$, $v=\Im w_n$, give local coordinates of $S^m$ with $e_0=0$.
Thus we may identify $J\calE$ with the space of real formal power series of $(w',\conj{w}',v)$. 
Each element of $J\calP$ is given by 
the real part of a formal power series of $\xi$ at $\wt e_0=(1,0,\dots,0)$ homogeneous 
of degree $0$.  Hence substituting $\xi_0=1,\xi'=w',\xi_n=|w'|^2/2+iv$,
we see that $J\calP$ consists of formal power series of the form
$$
\Re\sum_{|\alpha|\ge0,l\ge0}A_{\alpha}^l {w'}^\alpha (|w'|^2+2iv)^l,
\quad
\text{where } A_{\alpha}^l\in\C.
$$
To define a complementary space
of $J\calP$, we write $f\in J\calE$ as
\begin{equation}\label{A-exp}
f(w',\conj w',v)=\sum_{p,q\ge0}A_{p,q}(v),
\end{equation}
where
$$
A_{p,q}(v)=\sum_{
\scriptsize\begin{matrix}|\alpha|=p,|\beta|=q\\
l\ge0
\end{matrix}}A_{\alpha\conj\beta}^l
{w'}^\alpha {\conj{w}'}^\beta v^l.
$$
Let $\calN_0$ be the subspace (not a submodule) of $\calE$ 
defined by the equations
$$
A_{p,q}=0   \text{ if $\min(p,q)= 0$},
\quad \Delta' A_{1,1}=0, 
$$
where $\Delta'=\sum_{j=1}^{n-1}\pa_{w_j}\pa_{\conj{w}_j}$.
Then we have

\medskip

\begin{lemma} As a vector space
 $J\calE=\calN_0\oplus J\calP$.
\end{lemma}

\medskip

We have $J\calE/J\calP\cong\calN_0$ as vector spaces
and may define $P$-action on $\calN_0$ via this isomorphism.
We identify $\calN_0$ with the space of 
lists $(A_{\alpha\conj\beta}^l)$; so $P$ acts on the list
$(A_{\alpha\conj\beta}^l)$.
Then CR invariants
of weight $(-n,-n)$ are $P$-equivariant polynomial of
$(A_{\alpha\conj\beta}^l)$ that takes value in 
$\sigma_{-n,-n}$.
We define the weight of $A_{\alpha\conj\beta}^l$ to
be $|\alpha|+|\beta|+2l$ and extend the weight to 
the monomials of $(A_{\alpha\conj\beta}^l)$
by summing up the weight of each variable. 
It then turns out that  a CR invariant of weight $(-n,-n)$
is a polynomial of $(A_{\alpha\conj\beta}^l)$
homogeneous of weight $2n$ (examine the action of
dilation $(w',w_n)\mapsto(\lambda w', |\lambda|^2 w_n)$).
Since each nontrivial $A_{\alpha\conj\beta}^l$
has weight $\ge2$, we see that a nonlinear monomial
of homogeneous weight $2n$ depend only on the
variables
$(A_{\alpha\conj\beta}^l)$ of weight $\le2n-2$. 
Thus we may restrict our attention to the polynomials on
$$
\calN^{[2n]}_0=\{(A_{\alpha\conj\beta}^l)_{|\alpha|+|\beta|+2l<2n}:
(A_{\alpha\conj\beta}^l)\in\calN_0\}.
$$
Note that the kernel of the projection $\calN_0\to\calN^{[2n]}_0$
is a $P$-submodule and so $\calN^{[2n]}_0$ has a structure
of quotient $P$-module. As we already know the linear CR invariants,
we may reduce Problem 2 to 

\medskip

{\sc Problem 3.}
{\em  Write down all CR invariant $\calN^{[2n]}_0\to\sigma_{-n,-n}$. 
}
\medskip

In case $n=2$, we can determine CR invariants
of weight $(-n,-n)$ just by counting the weight of
the variables
$(A_{\alpha\conj\beta}^l)$. 
Observe that the condition
$\Delta' A_{1,1}=0$ implies $A_{1,1}=0$. Thus
$f\in\calN_0$ is of the form
$$
f=A_{2,1}+A_{1,2}+\sum_{p+q\ge4, p,q\ge1}A_{p,q}.
$$
It follows that each nontrivial $A_{\alpha\conj\beta}^l$
has weight $\ge3$. Thus nonlinear monomial should have
weight $\ge6$ and it cannot appear in an invariant of
weight $(-2,-2)$. Since
$A_{2,2}(0)$
is the only term of weight $4$, we see that
a CR invariant of weight $(-2,-2)$ must be
a multiple of $A_{2,2}(0)$.
(We know the existence of a CR invariant of
weight $(-2,-2)$; so $A_{2,2}(0)$
must be a CR invariant.) 

We next examine the
dependence of $\wt\nabla^{1,1}\wt f$
on $\calN_0$.
We use the coordinates $\xi_0,\dots,\xi_n$ to define the components
of 
$$\wt\nabla^{p,q}\wt f(\wt e_0)=(T_{I_1\cdots I_p\conj J_1\cdots \conj J_q }(\wt f)).
$$ 
The action of $P$ on $\wt f$ induces an action on
the tensors given by
$$
\mathbb{T}_0^{p,q}=\bigodot^p V^*\otimes\bigodot^q\conj V^*\otimes \sigma_{-p,-q},
$$
where $V=\C^{n+1}$ on which $P$ acts by right multiplication on column vectors.
We next introduce a weight on each component of $\mathbb{T}_0^{p,q}$.
Set $\|0\|=0$, $\|j\|=1$, $j=1,2,\dots, n-1$, $\|n\|=2$
and extend it to the list of indices by
$$\|I_1\cdots I_p \conj J_1 \cdots \conj J_q\|=
\|I_1\|+\cdots +\|I_p\| +\|J_1\|+\cdots +\|J_q\|.
$$
Since $\wt f$ is unique modulo $O(L^n)$, we have the
following estimate of the ambiguity. 

\medskip

\begin{lemma}
If $\|\calI \conj \calJ\|<2n$, then
$T_{\calI \conj \calJ}(\wt f)$
depends only on $f\in J\calE$ modulo
$J\calP$, where
$\calI=I_1\cdots I_p$ and $\conj\calJ=\conj J_1 \cdots \conj J_q$
with $p,q\ge1$.
\end{lemma}

\medskip

Let $\mathbb{T}_0=\prod_{p,q\ge1}\mathbb{T}_0^{p,q}$
and set 
$$
\mathbb{T}^{[2n]}_{0}=\{(T_{\calI \conj \calJ})_{
\|\calI \conj \calJ\|<2n}:
(T_{\calI \conj \calJ})\in\mathbb{T}_0\}.
$$
Since the kernel of the projection $\mathbb{T}_{0}\to\mathbb{T}^{[2n]}_{0}$ is shown to be a submodule, we may define $P$ action 
on $\mathbb{T}^{[2n]}_{0}$ as the quotient.
The lemma above ensures that the map
$$
T\colon J\calE/J\calP\to\mathbb{T}^{[2n]}_{0},
\quad T(f)=(T_{\calI \conj \calJ}(f))_{
\|\calI \conj \calJ\|<2n},
$$
is a well-defined $P$-equivariant map.
It is easy to verify that the image of $T$ satisfies the 
equations
\begin{equation}\label{harmonic-eq}
\conj{T_{\calI \conj \calJ}}=T_{\calJ\conj \calI },\quad
{\wt g}^{K\conj L}
T_{\calI K\conj L\,\conj\calJ }=0,\quad
T_{\calI 0\conj\calJ }=-
|\calI|T_{\calI \conj\calJ },
\end{equation}
where $|I_1\cdots I_p|=p$ and we also consider the case $p=0$. 
These respectively follow from the facts that $\wt f$ is
real, $\Delta$-harmonic and homogenous of degree $(0,0)$.
Let us denote by $\calH_0$ the submodule of $\mathbb{T}_0$ defined
by these equations and define $\calH_0^{[2n]}$ to be its
projection to $\mathbb{T}^{[2n]}_{0}$.
It is clear that $T(f)$ depends only on finite jets of
$J\calE/J\calP\cong\calN_0$. 
More precisely, we can show that $T_{\calI \conj \calJ}(f)$
of weight $w$ is a homogeneous polynomial of
$(A_{\alpha\conj\beta}^l)$ of weight $w$. It follows that
$T\colon\calN^{[2n]}_0\to\calH^{[2n]}_0$ is well-defined.

\medskip

\begin{theorem}[Jet Isomorphism Theorem]\label{JIT}
The map $T\colon\calN^{[2n]}_0\to\calH^{[2n]}_0$
is an isomorphism of $P$-modules.
\end{theorem}

\medskip

The isomorphism $T$ reduces the problem of determining 
invariants on $\calN^{[2n]}_0$ to the one on 
$\calH^{[2n]}_0$.  It is clear that
an invariant polynomial on $\calH^{[2n]}_0$
can be also seen as an invariant polynomial
on $\calH_0$.  The converse is also true
because an invariant taking values in 
$\sigma_{-n,-n}$ depends only on the components
in $\calH^{[2n]}_0$. Hence Problem 3 is reduced to
\medskip

{\sc Problem 4.}
{\em  Write down all CR invariant
$\calH_0\to\sigma_{-n,-n}$. 
}
\medskip

Fortunately, this problem has been completely solved in more general
setting.

\medskip

\begin{theorem}[{[BEG]}] \label{BEGthm}Let 
$I\colon\calH_0\to\sigma_{q,q}$ be a CR invariant for some integer $q$. Then $I$ is a linear combination
of complete contractions of the form
$$
 \contr(T^{p_1,q_1}\otimes \cdots\otimes T^{p_k,q_k}),
 \quad T^{p,q}\in \mathbb{T}_0^{p,q}.
$$
\end{theorem}

Note that $T^{p,q}$ is trace-free and there is no linear CR invariant
defined on $\calH_0$.
On the other hand, the linear CR invariant operator
$Q(f)$ is the complete contraction of $T^{n,n}(\wt f)$, which is not covered by the jet isomorphism theorem.
Thus there is no hope to generalize Theorem \ref{JIT} 
to $T\colon\calN^{[w]}_0\to\calH^{[w]}_0$, $w>2n$.
To get some relation between 
$\calN_0$ and $\calH_0$ for higher jet, one needs to introduce
additional parameter to resolve the ambiguity of the
harmonic extension. This approach is discussed in \cite{H2}.

\section{Bernstein-Gelfand-Gelfand resolution}

We explain why it is natural to realize the quotient space
$J\calE/J\calP$ as tensor space in $\T_0$ 
from the point of view of the BGG resolution.
We here only consider some explicit examples of (the dual of) BGG resolutions; more complete treatment of BGG should be contained in other articles in this volume.
We first assume that $n>2$. The case $n=2$ will be discussed at the end of this section.

Recall that the complexified tangent space of $S^m$ admits 
a natural subbundle $T^{1,0}_b=\C T S^m\cap T^{1,0}\C^n$ and its conjugate $T^{0,1}_b$. Hence the restrictions of the exterior derivative of a function give $\pa_b f=df|_{T^{1,0}}$ and
$\conj \pa_b f=df|_{T^{0,1}}$.  Let $\calE^{1,0}$ and $\calE^{0,1}$ denote the space of the sections of $(T^{1,0}_b)^*$ and $(T^{0,1}_b)^*$ respectively. For $p+q\le n-1$, we also define
$\calE^{p,q}$ to be the space of the sections of $\bigwedge^p (T^{1,0}_b)^*\otimes_\circ\bigwedge^q (T^{0,1}_b)^*$,
where $\otimes_\circ$ means the trace-free part of the bundle
with respect to the Levi from. Note that these bundles are induced from irreducible representations of $P$.  We fix a subbundle $N$ such that
$$
\C TS^m=N\oplus T^{1,0}\oplus T^{0,1}
$$
and identify $f\in\calE^{p,q}$ with a $(p+q)$-form. 
If $p+q<n$, we can define
$$
\calE^{p,q}
\begin{matrix}
\overset{\pa_b }{
\nelongarrow
}
  \\ 
\underset{\conj\pa_b }{
\selongarrow
} 
\end{matrix}
\quad
\begin{matrix}
\calE^{p+1,q}\\ \  
\\
\calE^{p,q+1}
\end{matrix}
$$
by applying exterior derivative to the $(p+q)$-form and projecting to appropriate bundles; these maps are independent of the choice of $N$.

The decompositions of the exterior derivatives on functions and one-forms
give a complex of $G$-modules ($\C$ is regarded as the trivial representation)
\begin{equation}\label{KR}
0\to\mathbb{C}\to
\calE^\bC
\begin{matrix}
\nelongarrow
  \\ 
\selongarrow 
\end{matrix}
\quad
\begin{matrix}
\calE^{0,1}\\ \  
\\
\calE^{1,0}
\end{matrix}
\quad
\begin{matrix}
\nelongarrow
  \\ 
\selongarrow
\\
\nelongarrow
  \\ 
\selongarrow
\end{matrix}
\qquad
\begin{matrix}
\calE^{0,2} \\ \
\\
\calE^{1,1} 
\\ \ \\
\calE^{2,0} 
\end{matrix}
\end{equation}
Here $\calE^\C$ is the space of complex valued smooth functions.
By composition, we can also define a $G$-equivariant map
$$
R\colon\calE^\C\to\calE^{1,1},
\quad R(f)=\text{trace-free part of }i\,\conj\pa_b\pa_bf.
$$ 
Restricting $R$ to $\calE$ and
then taking the jets at $e_0$, we obtain a homomorphism of $P$-modules
$$
R:J\calE\to J\calE^{1,1}.
$$
The kernel is given by $J\calP$ (see e.g. \cite{L}) and we get an injection
$$
J\calE/J\calP \to J\calE^{1,1}.
$$
The map $T\colon J\calE/J\calP\to\T_0^{[n]}$ is given by
taking the jet at $\wt e_0$ of $\wt\nabla^{1,1}\wt f$ and
$\wt\nabla^{1,1}$ can be seen as 
the lift of $R$ to the ambient space. 
In terms of the slice $\calN_0$, the value of $R(f)$ at $e_0$
is given by $A_{1,1}(0)$, the leading term of the normal form.
Since $\Delta' A_{1,1}=0$, it defines a trace-free tensor.

We next realize $J\calE/J\calP$ as
a cohomology of a subcomplex of \eqref{KR}. 
Since
$J\calE^\C\to J\calE^{0,1}\to J\calE^{0,2}
$ is exact (see e.g.\ \cite{Bo}), it gives no local invariants. However, 
by replacing $\calE^\C$ by $\calE$, we obtain a complex
$$
0\to\bR\to J\calE\to J\calE^{0,1}\to J\calE^{0,2}
$$
which has nontrivial first cohomology
$$
H^1=\frac{\ker (\conj\pa_b:J\calE^{0,1}\to J\calE^{0,2})}{
\image (\conj\pa_b:J\calE\to J\calE^{0,1})}.
$$
We will show that this is isomorphic to $J\calE/J\calP$.

\medskip

\begin{proposition}\label{H1-isom}  The map $J\calE\to J\calE^{0,1}$ induces
an isomorphism
$
J\calE/J\calP\cong H^1
$.
\end{proposition}

\smallskip

\begin{proof}
By the exactness of $J\calE^\bC\to J\calE^{0,1}\to J\calE^{0,2}$,
we may find for each $f\in\ker (\conj\pa_b\colon J\calE^{0,1}\to J\calE^{0,2})$
a jet of function $u+iv\in J\calE^\bC$ such that $f=\conj\pa_b(u+iv)\equiv i\conj\pa_bv$
mod $\conj\pa_b(J\calE)$. 
Thus $i\conj\pa_b\colon \calE\to H^1$ is surjective. It suffices to compute the kernel. If $i\conj{\pa}_b v\in\conj{\pa}_bJ\calE$,
we may take $u\in J\calE$ so that
 $i\conj{\pa}_b v=\conj{\pa}_b u$, or equivalently, $\conj{\pa}_b (u-iv)=0$. Thus $v\in J\calP$.
 \end{proof}

\medskip

We next apply the construction of the quotient module to the BGG resolution of the adjoint representation
$\mathfrak{su}(n,1)$.
We first complexify the representation and give the
BGG for $\mathfrak{sl}(n+1,\bC)$:
\begin{equation*}
0\to\mathfrak{sl}(n+1,\bC)\to
\calE^\bC(1)
\begin{matrix}
%\overset{D_{0,0}^+\ }{
\nelongarrow
%}
  \\ 
%\underset{D_{0,0}^-\ }{
\selongarrow
%} 
\end{matrix}
\quad
\begin{matrix}
\calE^{(1,1)}(1)\\ \  
\\
\conj{\calE^{(1,1)}(1)}
\end{matrix}
\quad
\begin{matrix}
%\overset{D_{1,0}^+\ }
{\nelongarrow}
  \\ 
%\underset{\scriptscriptsize{D_{1,0}^-}\ }
{\selongarrow} 
\\
%\overset{D_{0,1}^+\ }
{\nelongarrow}
  \\ 
%\underset{D_{0,1}^-\ }
{\selongarrow} 
\end{matrix}
\qquad
\begin{matrix}
\calE^{(2,1)}(1)
\\ \
\\
\calE^{(1,1)(1,1)}(1)
\\ \ \\
\conj{\calE^{(2,1)}(1)}
\end{matrix}
\end{equation*}
Here $\calE^{(p,q)}(1)$ denotes the space of the sections of
$\bigwedge^p (T^{0,1})^*
\otimes\bigwedge^q (T^{0,1})^*\otimes\sigma_{1,1}$
satisfying
$
 f_{[\alpha_1\cdots\alpha_p\beta_1]\beta_2\cdots\beta_q}=0,
$
namely, $f_{\alpha_1\dots\beta_q}$ has the symmetry corresponds to the Young diagram with two
columns of hight $p$ and $q$.
$\calE^{(1,1)(1,1)}(1)$ is the space of the sections of
 $\bigodot^2 (T^{1,0})^*
\otimes_\circ\bigodot^2 (T^{0,1})^*\otimes\sigma_{1,1}$.

Composing the homomorphisms in the complex, 
we obtain a map
$$
 R\colon\calE^\bC(1)\to\calE^{(1,1)(1,1)}(1).
$$
It has order 4 and is given by the projection to 
$\calE^{(1,1)(1,1)}(1)$ of
$\nabla^4 f$, where $\nabla$ is the Tanaka-Webster connection.
The kernel of the map
$$
R\colon J\calE(1)\to J\calE^{(1,1)(1,1)},
$$
denoted by $J\calP_1$,
is given by the jets of the form
$$
 \Re\sum_{j=0}^n\conj\xi_j f_j(\xi),
$$
where $f_j$ is a formal power series of $\xi$ at $\wt e_0$ of homogeneous
degree $1$.
In this setting, we can say that 
the formal theory of Moser's normal form gives a slice
of the coset of $J\calE(1)/J\calP_1$.
Moser's normal form
is a slice of equivalence classes of
the real hypersurfaces of the form
$$
2\Re w_n=|w'|^2+f(w',\conj{w}',v)
$$
modulo holomorphic coordinates changes that fix $0$.
$J\calP_1$ describes the linear part of the effect by
coordinates changes. By setting $\xi_0=1$, we identify
$J\calE(1)$ with the space of real formal 
power series $f(w',\conj{w}',v )$ of the form
\eqref{A-exp}.  Then the space of Moser's normal form
$\calN_1$ is given by the conditions:
$$
\begin{aligned}
 A_{p,q}&=0 \quad \text{ if $\min(p,q)\le 1$},
 \\
 (\Delta')^{j+k+1} A_{2+j,2+k}&=0
 \quad \text{for $0\le j,k\le 1$}.
\end{aligned}
$$

\begin{proposition}[{\cite{CM}}]\label{Moser}
As a vector space $J\calE(1)=\calN_{1}\oplus J\calP_1$.
\end{proposition}

\medskip

If we identify $J\calE(1)/ J\calP_1\cong\calN_{1}$,
the value of $R(f)$ at $e_0$ is given by $A_{2,2}(0)$,
which can be idenfitied with a trace-free tensor because of the normalization condition
$\Delta'A_{2,2}(0)=0$.

We next realize $J\calE(1)/J\calP_1$
as a cohomology.
If we start with the adjoint representation,
we obtain
$$
0\to\mathfrak{su}(1,n)\to\calE(1)
\to
\calE^{(1,1)}(1)
\quad
\begin{matrix}
{\nelongarrow}
  \\ 
{\selongarrow} 
\end{matrix}
\qquad
\begin{matrix}
\calE^{(2,1)}(1)
\\ \
\\
H\calE^{(1,1)(1,1)}_\circ(1)
\end{matrix}
$$
where $H\calE^{(1,1)(1,1)}_\circ(1)$ is the subspace
of $\calE^{(1,1)(1,1)}_\circ(1)$ consisting of tensors with
hermitian symmetry.
The subcomplex 
$$
0\to\calE(1)\overset{D_0}{\longrightarrow}\calE^{(1,1)}(1)
\overset{D_1}{\longrightarrow}
\calE^{(2,1)}(1)\to\cdots\to\calE^{(n-1,1)}(1)\to 0
$$
is called the deformation complex (see\cite{Cap}): $\ker D_1$
parametrizes infinitesimal integrable deformation
of CR structure of $S^m$, and $\image D_0$ corresponds to
trivial deformations.
In analogy with Proposition \ref{H1-isom} we have
$$
J\calE(1)/J\calP_1\cong \frac{\ker(D_1:J\calE^{(1,1)}(1)\to J\calE^{(2,1)}(1))}
{\image 
(D_0:J\calE(1)\to J\calE^{(1,1)}(1))}.
$$

We finally consider the case $n=2$.
For the trivial representation, we have the following complex
of CR invariant differential operators
$$
0\to\C\to\calE^\C\to\begin{matrix}
\calE^{0,1}\\ \oplus
\\
\calE^{1,0}
\end{matrix}
\to\begin{matrix}
\calE^{0,1}(-1)\\ \oplus
\\
\calE^{1,0}(-1)
\end{matrix}
\to\calE^\C(-2)\to0,
$$
which is called the Rumin complex.
The operator in the middle has order $2$
and the composition
$$ 
R\colon\calE\to\calE^{0,1}(-1)
$$
has order 3.  
$R(f)=0$ characterizes CR pluriharmonic functions
and $J\calE/J\calP$ can be embedded into $\calE^{0,1}(-1)$.
Under the identification $J\calE/J\calP\cong \calN_0$,
the value of $R(f)$ at $e_0$ is given by $A_{1,2}(0)$,
which is the third order term of $f$.
It seems rather strange that the one form valued operator $R$
lifted to the ambient space as a $(1,1)$-tensor valued operator.
This type of change also appears in the ambient metric construction of 3-dimensional conformal structures,
where Cotton tensor lifted to the ambient curvature tensor.

For the adjoint representation, we have
$$
0\to\mathfrak{sl}(3,\bC)\to\calE^\C(1)\to
\begin{matrix}
\calE^{(1,1)}(1)
\\ \oplus
\\
\conj{\calE^{(1,1)}(1)}
\end{matrix}
\to
\begin{matrix}
\calE^{(1,1)}(-1)\\ \oplus
\\
\conj{\calE^{(1,1)}(-1)}
\end{matrix}
\to\calE^\C(-3)\to0.
$$
The operator in the middle has order 4 and the composition $\calE(1)\to\conj{\calE^{(1,1)}(-1)}$ has order $6$.
Under the identification $J\calE(1)/J\calP_1\cong \calN_1$,
the value of $R(f)$ at $e_0$ is given by $A_{2,4}(0)$.
This is the leading term
of Moser's normal form and has order $6$.
The ambient lift of $R$ is $\wt\nabla^{2,2}$;
it again change the type of the tensor.
 
 \section{Concluding remarks}
 The arguments in section 5 can be applied to
 CR invariants $J\calE/J\calP\to\sigma_{-p,-p}$
 for $p=0,1,2,\dots, n$. These invariants are given by complete
 contractions on the ambient space. 
 To study the case $p>n$, we need to generalize the
 jet isomorphism theorem. One way of doing this is
 to allow logarithmic term in the solution 
 to $\Delta \wt f=0$. Since the solution is not unique, we need to
 introduce a parameter space to formulate the jet isomorphism
 theorem. This causes a new difficulty in the construction of
 CR invariants and the classification problem in this case is still open \cite{H2}. 
 
 To deal with the general strictly pseudoconvex real hypersurfaces
 in a complex manifold $X$, we need to start with
 the construction of the ambient space.  The ambient metric
 is defined on the canonical bundle of $X$ (with the zero section removed) and is constructed by solving a complex \MA equation. The ambient metric construction in Section 4 also holds for this curved ambient metric. The details will be given in my forthcoming paper.

\end{document}